\documentclass[12pt]{article}
\usepackage[a3paper]{geometry}
\usepackage[numbers]{natbib}
\usepackage{graphicx}
\usepackage{color, multirow,amsfonts,amsmath}
\usepackage{url}
\usepackage{rotating}
\begin{document}
\title{Merging of bivariate \\compound Binomial processes with shocks\\({\it{Working version}})
\\ Pavlina K. Jordanova
\\{\small{\it Faculty of Mathematics and Informatics, Konstantin Preslavsky University of Shumen, \\115 "Universitetska" str., 9712 Shumen, Bulgaria.
\\ Corresponding author: pavlina\_kj@abv.bg}}
\\ Evelina Veleva
\\{\small{\it{Department of Applied mathematics and Statistics, "Angel Kanchev" University of Ruse, Bulgaria.}}}}

\maketitle

\begin{abstract}
 The paper investigates a discrete time Binomial risk model with different types of polices and shock events may influence some of the claim sizes. It is shown that this model can be considered as a particular case of the classical compound Binomial model. As far as we work with parallel Binomial counting processes in infinite time, if we consider them as independent, the probability of the event they to have at least once simultaneous jumps would be equal to one. We overcome this problem by using thinning instead of convolution operation. The bivariate claim counting processes are expressed in two different ways. The characteristics of the total claim amount processes are derived. The risk reserve process and the probabilities of ruin are discussed. The deficit at ruin is thoroughly investigated when the initial capital is zero. Its mean, probability mass function and probability generating function are obtained. We show that although the probability generating function of the global maxima of the random walk is uniquely determined via its probability mass function and vice versa, any compound geometric distribution with non-negative summands has uncountably many stochastically equivalent compound geometric presentations. The probability to survive in much more general settings,  than those, discussed here, for example in the Anderson risk model, has uncountably many Beekman's convolution series presentations.
\end{abstract}

\section{INTRODUCTION}

The discrete time classical Binomial risk model is considered for example in H. Gerber in 1988 \cite{Gerber1988}. The total claim amount process is a compound Binomial processes. The claim sizes are discrete, and the premium income rate $c = 1$. The investigation of the corresponding probabilities of ruin and non-ruin originates since his works and those of E.S.W. Shiu (1989) \cite{Shiu1989}. The investigation of probabilities of ruin in these models can be reduced to the one of the distribution of the global maxima of the corresponding random walk. Therefore, they have Beekman's (1968) \cite{Beekman1968} convolution series (compound geometric) presentation. A very good explanation of risk theory results as a particular case of  the random walk theory could be seen in Rolski et al. (1998) \cite{Rolski1998} or Feller (1957,2008) \cite{Feller1957,Feller2008}.

The model considered in this work includes possibility to have two different types of polices and shock events to influence both of them. It is shown that the considered model can be reduced to a compound Binomial risk model.

According to the zero-ones law, if we observe some number of independent binomial processes in infinite time, the probability  of the event they to have at least once simultaneous jumps is equal to one. We overcome this problem by using thinning instead of convolution operation. The bivariate claim counting processes are expressed in two different ways. Their probability generating functions (p.g.fs.) show that they have bivariate Binomial distribution. The last is considered for example in Johnson et al. (1997) \cite{Johnsonetal1997}. The Laplace-Stieltjes transforms (LSTs) of the total claim amount processes and the corresponding  numerical characteristic are obtained. Finally the probabilities of ruin and distribution of the deficit at ruin are derived.  In an analogous way the results can be obtained for more types of polices or different types of shock events.

We assume that all random elements discussed here are measurable with respect to one the same probability space  $(\Omega, \mathcal{F}, \mathbb{P})$ with natural filtration born by the considered processes. Along the paper consecutively applied symbol $\bot$  means mutually independent. The letter $\mathbb{G}$ is for probability generating function (p.g.f.). For example the p.g.f. of the random vector $(\xi_1, \xi_2)$ is denoted by $\mathbb{G}_{\xi_1, \xi_2}(z_1, z_2) = \mathbb{E}(z_1^{\xi_1}z_2^{\xi_2})$. In the univariate case, it is the p.g.f. of the corresponding random variable (r.v.). Analogously, the letter $\mathbb{L}$ means LST. For any $(\Omega, \mathcal{F}, \mathbb{P})$  measurable event $A$, $I_{\{A\}}$ denotes a Bernoulli r.v. with parameter $\mathbb{P}(A)$.

\section{DESCRIPTION OF THE MODEL}

Let $B_0, B_1, B_2 = \{(B_0(t), B_1(t), B_2(t)), \, t  = 0, 1, ...\}$ be a Multimonial process, with independent and homogeneous additive increments and parameters correspondingly $p_0, p_1, p_2 \geq 0$, $p: = p_0 + p_1 + p_2 \leq 1$. This process is determined via the distribution of its time intersections
\begin{equation}\label{multinomial}
\mathbb{P}(B_0(t) = b_0, B_1(t) = b_1, B_2(t) = b_2) = \frac{t!}{b_0!b_1!b_2!(t - b_0 - b_1 - b_2)!}p_0^{b_0}p_1^{b_1}p_2^{b_2}(1-p)^{t - b_0 - b_1 - b_2},
\end{equation}
for $b_0 + b_1 + b_2 \leq t$, and $b_0, b_1, b_2 \in \{0, 1, ..., t\}$, and this probability mass function (p.m.f.) is equal to $0$, otherwise. The  properties of Multinomial distribution could be found in many textbooks in probability, for example, N. L. Johnson et al. (1997) \cite{Johnsonetal1997}.

We suppose that the claims can be of two mutually exclusive different types. There is a possibility for covering of shock events (like car-crashes) which cause payment of claims in both polices. The counting process of the shock events is denoted by $B_0 = \{B_0(t), \, t  = 0, 1, ...\}$.
For any fixed $k = 1, 2, 3, 4$ we suppose that $Y_{k1}, Y_{k2}, ...$ are independent integer-valued identically distributed (i.i.d.) random variables (r.vs.) with cumulative distribution function (c.d.f.) $\mathbb{F}_k$ and $\{Y_{1i}\}_{i=1}^\infty \bot \{Y_{2i}\}_{i=1}^\infty \bot \{(Y_{3i},Y_{4i})\}_{i=1}^\infty$. The r.v. $Y_{ki}$ models the $i$-th claim amount of a certain type $k$. Without lost of generality $\mathbb{P}(Y_{ki} > 0) = 1$, otherwise we will change the parameter of the corresponding counting binomial process, participating in $S_1$ and $S_2$, defined in (\ref{S1S2}). For $k = 1$, $Y_{1i}$ is the part which is due to $i$-th standard insurance policy of type $1$ exclusively the common shocks. Analogously for $k = 2$. $Y_{3i} + Y_{4i}$ is the part which have to be paid to the customer due to some common shock. For example, if we suppose that $k = 1$ is for  health insurance, and $k = 2$ is for car insurance, then $\sum_{i = 0}^{B_1(t)} Y_{1i}$ is the total claim amount up to time $t$, which is due to some healthy problems, not related with car-crashes. $\sum_{i = 0}^{B_0(t)} Y_{3i}$ is the total claim amount up to time $t$, which is due to healthy problems, caused by car-crashes. $\sum_{i = 0}^{B_2(t)} Y_{2i}$ is the total claim amount up to time $t$, which is due to some car problems which are not related with car-crashes, for example  car-thefts.
$\sum_{i = 1}^{B_0(t)} Y_{4i}$ is the total claim amount up to time $t$, which is due to some car problems which are related with car-crashes.

In general, the  counting processes $M_1$ and $M_2$ of both types of claims are dependent Binomial process with common shocks. More precisely for all $ t \geq 0$
\begin{equation}\label{M}
M_1(t) = B_1(t) +  B_0(t), \quad M_2(t) = B_2(t) +  B_0(t),
\end{equation}
 We are going to describe them in the next section. Note that for $p_0 \in (0, 1)$, the processes $B_1, B_2$ and $B_0$ are always dependent.

$S = \{S(t):  t = 0, 1, ...\}$ is the total claim amount process.   For all $t = 0, 1, ...$, it satisfies the following relations
 \begin{eqnarray}
\label{S}   S(t) &=& S_1(t) + S_2(t),\\
\label{S1S2} S_1(t) &=& \sum_{i = 0}^{B_1(t)} Y_{1i} + \sum_{i = 0}^{B_0(t)} Y_{3i},\quad Y_{10} = 0,\,\, Y_{30} = 0\\
\nonumber S_2(t) &=& \sum_{i = 0}^{B_2(t)} Y_{2i} + \sum_{i = 0}^{B_0(t)} Y_{4i},\quad Y_{20} = 0,\,\, Y_{40} = 0.\\
\nonumber \{Y_{ki}\}_{i = 1}^\infty &\bot& \{B_r(s), s\geq 0\},  \,\, k = 1, 2, 3, 4, \,\, r = 0, 1, 2.
 \end{eqnarray}
The processes $S_1 = \{S_1(t), \,  t = 0, 1, ...\}$ and $S_2 = \{S_2(t), \, t = 0, 1, ...\}$ are dependent. In our example $S_1$ is the total claim amount processes for the events caused by "healthy problems", $S_2$ is those of the events caused by "problems with cars".

Finally we will consider a risk reserve process $R_u = \{R_u(t): t \geq 0\}$ defined via the equality
\begin{equation}\label{R}
R_u(t) = u + t - S(t), \quad t = 0, 1, ...,
\end{equation}
where $u = 0, 1, ...$ is the initial capital. The corresponding time of ruin is according to the definition of Shiu (1989) \cite{Shiu1989}
\begin{eqnarray}
\label{timeofruin}
\tau(u)&:=& \inf\{t > 0: R_u(t) < 0\},\quad  \inf\, \emptyset = \infty,
\end{eqnarray}
and the probability for ruin in "infinite horizon" and initial capital $u$ as usually will be denoted by
\begin{equation}\label{probforruin}
\psi(u) = \mathbb{P}(\tau(u) < \infty).
\end{equation}
The corresponding survival probability is $\delta(u) = 1 - \psi(u)$ .

\bigskip
\section{THE BIVARIATE COUNTING PROCESS}

Binomial and bivariate compound Binomial processes are very well investigated in the scientific literature. Let us now derive compound Binomial presentation of the bivariate counting process $(M_1, M_2) := \{(M_1(t), M_2(t)): t \geq 0\}$.

In order to formulate our results let us define a random vector $(I_1, I_2, I_0)$ with probability mass function (p.m.f.)
\begin{equation}\label{pmf1}
\mathbb{P}(I_1 = 1, I_2 = 0, I_0 = 0) = \frac{p_1}{p}, \,\,  \mathbb{P}(I_1 = 0, I_2 = 1, I_0 = 0) = \frac{p_2}{p}, \,\,\mathbb{P}(I_1 = 0, I_2 = 0, I_0 = 1) = \frac{p_0}{p},
\end{equation}
and zero otherwise. Assume that the random vectors $(I_{11}, I_{21}, I_{01}), (I_{12}, I_{22}, I_{02}),  ...$ are i.i.d. with p.m.f. (\ref{pmf1}).

Denote by, $B = \{B(t), \, t = 0, 1, ...\}$ a homogeneous Binomial process with independent additive increments and parameter $p$, and by
\begin{eqnarray}
A_1(t) &:=& \sum_{i = 0}^{B(t)} (I_{1i} + I_{0i}) = \sum_{i = 0}^{t} (I_{A_i}(I_{1i} + I_{0i})),\label{A11}\\
A_2(t) &:=& \sum_{i = 0}^{B(t)} (I_{2i} + I_{0i}) = \sum_{i = 0}^{t} (I_{A_i}(I_{2i} + I_{0i})), \quad t = 0, 1, ....\label{A21}
\end{eqnarray}
Here $I_{10} = I_{20} = I_{00} = 0$,  $B$ and $(I_{11}, I_{21}, I_{01}), (I_{12}, I_{22}, I_{02}),  ...$ are independent, for any fixed $i = 1, 2, ...$, $\mathbb{P}(A_i) = p$, and $I_{A_i}$ are independent on $(I_{1i}, I_{2i}, I_{0i})$. The event $A_i$ can be interpreted as "in the time point $t = i$ there is an insurance event". In our example the random process $A_1$ counts the number of the insurance events caused by healthy problems. The random process $A_2$ counts the number of the insurance events caused by problems with cars.

In the next theorem we show that the bivariate counting process
$(M_1, M_2) := \{(M_1(t), M_2(t)): t = 0, 1, ...\}$ is a particular case of multivariate compound Binomial processes with equal number of summands, i.e. of type I, in the sense of Sundt and Vernic \cite{SundtVernic, Vernic2018}.

{\bf Theorem 1.} The bivariate processes $(M_1, M_2)$ and $(A_1, A_2)$ coincide in the sense of their finite dimensional distributions.

{\bf Proof:} The processes $(M_1, M_2)$ and $(A_1, A_2)$ have homogeneous and independent additive increments and they start from the coordinate beginning. Therefore, in order to prove their stochastic equivalence it is enough to prove equality of their univariate time intersections. Due to the uniqueness of the correspondence between the probability laws and their p.g.fs. it is enough to derive equality between the p.g.fs. The definition of p.g.fs., (\ref{M}) and the definition for  Multinomial distribution of $(B_0(t), B_1(t), B_2(t))$ entail that for all $t \in \mathbb{N}$,
\begin{eqnarray*}
&\mathbb{E}&\left[z_1^{M_1(t)}z_2^{M_2(t)}\right]\\
&=&\mathbb{E}\left[z_1^{B_1(t) +  B_0(t)}z_2^{B_2(t) +  B_0(t)}\right]  = \mathbb{E}\left[z_1^{B_1(t)}z_2^{B_2(t)} (z_1z_2)^{B_0(t)} \right]  \\
&=&  \sum_{\substack{n_0, n_1, n_2 \in \{0, 1, ..., t\},\\ n_1 + n_2  + n_0 \leq t}} z_1^{n_1} z_2^{n_2} (z_1z_2)^{n_0}\mathbb{P}(B_1(t)=n_1,B_2(t)=n_2,B_0(t)=n_0)  \\
&=&  \sum_{\substack{n_0, n_1, n_2 \in \{0, 1, ..., t\},\\ n_1 + n_2  + n_0 \leq t}} z_1^{n_1} z_2^{n_2}(z_1z_2)^{n_0}\frac{t!p_1^{n_1}p_2^{n_2}p_0^{n_0}(1 - p_1 - p_2 - p_0)^{t - n_1 -n_2 - n_0}}{n_1!n_2!n_0!(t - n_1 - n_2 - n_0)!}  \\
&=& (1 - p)^{t} \sum_{\substack{n_0, n_1, n_2 \in \{0, 1, ..., t\},\\ n_1 + n_2  + n_0 \leq t}} \left[\frac{p_1}{1-p}z_1\right]^{n_1} \left[\frac{p_2}{1-p}z_2\right]^{n_2}\left[\frac{p_0}{1-p}(z_1z_2)\right]^{n_0}\frac{t!}{n_1!n_2!n_0!(t - n_1 - n_2 - n_0)!}\\
&=& (1-p)^t \left[1 + \frac{p_1}{1-p}z_1 + \frac{p_2}{1-p}z_2 + \frac{p_0}{1-p}z_1z_2\right]^t\\
&=& \left[1 - p + p_1z_1 + p_2z_2 + p_0z_1z_2\right]^t.
\end{eqnarray*}

By the definitions (\ref{A11}), (\ref{A21}) and the formula for double expectations
\begin{eqnarray*}
\mathbb{E}\left[z_1^{A_1(t)}z_2^{A_2(t)}\right] &=&  \mathbb{E}\left[z_1^{\sum_{i = 0}^{B(t)} (I_{1i} + I_{0i})} z_2^{\sum_{i = 0}^{B(t)} (I_{2i} + I_{0i})} \right]\\
&=& \sum_{n=0}^\infty \mathbb{E}\left[z_1^{\sum_{i = 0}^{n} (I_{1i} + I_{0i})} z_2^{\sum_{i = 0}^{n} (I_{2i} + I_{0i})} \right]\mathbb{P}(B(t) = n)
\end{eqnarray*}

By the multiplicative property of p.g.fs. and the fact that the random vectors $(I_{11}, I_{21}, I_{01}), (I_{12}, I_{22}, I_{02}),  ...$ are i.i.d. we have
\begin{equation}
\mathbb{E}\left[z_1^{A_1(t)}z_2^{A_2(t)}\right] = \sum_{n=0}^\infty \left[\mathbb{E}\left(z_1^{I_{1i} + I_{01}} z_2^{I_{21} + I_{01}} \right)\right]^n \mathbb{P}(B(t) = n).
\end{equation}

Now, (\ref{pmf1}) entails
\begin{equation}
\mathbb{E}\left[z_1^{I_{1i} + I_{01}} z_2^{I_{21} + I_{01}} \right] = \frac{p_1}{p}z_1 + \frac{p_2}{p}z_2 + \frac{p_0}{p}z_1z_2.
\end{equation}

Therefore, the fact that $B$ is a Binomial process with parameter $p$ leads us to
\begin{eqnarray*}
\mathbb{E}\left[z_1^{A_1(t)}z_2^{A_2(t)}\right]&=& \sum_{n=0}^\infty \left[\frac{p_1}{p}z_1 + \frac{p_2}{p}z_2 + \frac{p_0}{p} z_1z_2\right]^n \mathbb{P}(B(t) = n)\\
&=&  \left[1 - p + p\left(\frac{p_1}{p}z_1 + \frac{p_2}{p}z_2 + \frac{p_0}{p}z_1z_2\right)\right]^t\\
&=&  \left[1 - p + p_1z_1 + p_2z_2 + p_0z_1z_2\right]^t
\end{eqnarray*}
which is exactly $\mathbb{E}\left[z_1^{M_1(t)}z_2^{M_2(t)}\right]$. The uniqueness of the correspondence between m.g.fs. and the probability distribution completes the proof. \hfill $\Box$

{\it Note 1.} For all $t = 0, 1, ...$ the distribution of the time intersections $(M_1(t), M_2(t))$ of the bivariate counting process $(M_1, M_2)$, defined in (\ref{M}) are determined via the following p.g.fs.,
$$\mathbb{G}_{M_1(t),M_2(t)}(z_1, z_2) = \mathbb{E}(z_1^{M_1(t)}z_2^{M_2(t)}) = \left[1 - p + p_1z_1 + p_2z_2 + p_0z_1z_2\right]^t.$$
These distributions are called bivariate Bernoulli distribution and they are very well investigated in the scientific literature. See for example M. A. Hamdan and H.A.Al-Bayyati (1969) \cite{HamdanAlBayyati1969} and Johnson et al. (1997) \cite{Johnsonetal1997}. For $m_1 , m_2 \in \{0, 1, ..., t\}$, and $m_1 + m_2 \leq t$
$$\mathbb{P}(M_1(t) = m_1, M_2(t) = m_2) = \sum_{i = 0}^{min(m_1, m_2)} \frac{t!p_0^ip_1^{m_1 - i}p_2^{m_2 - i}(1-p)^{t - m_1 - m_2 + i}}{i!(m_1 - i)!(m_2 - i)!(t - m_1 - m_2 + i)!},$$
and this probability mass function (p.m.f.) is equal to $0$ otherwise.

Thus, for $(r,s) = (1,2)$ or $(r,s) = (2,1)$, $m_s \in [0, t] \cup \mathbb{N}$ and $m_r \in [0, t-m_s] \cup \mathbb{N}$
$$\mathbb{P}(M_r(t) = m_r|M_s(t) = m_s)= \sum_{i = 0}^{min(m_r, m_s)} \left(\begin{array}{c}
                                                                              m_s \\
                                                                              i
                                                                            \end{array}
\right)\left(\begin{array}{c}
                                                                              t-m_s \\
                                                                              m_r-i
                                                                            \end{array}
\right) \left[\frac{p_0(1-p)}{p_rp_s}\right]^i\left(\frac{p_r}{1 - p}\right)^{m_r}\left[\frac{p_s(1-p_0-p_s)}{(p_0+p_s)(1 - p)}\right]^{m_s}\left(\frac{1-p}{1 - p+p_r}\right)^{t},$$
and this probability mass function (p.m.f.) is equal to $0$ otherwise.

The corresponding  mean square regression is linear. More precisely,
$$\mathbb{E}(M_s(t)|M_r(t)=i)= \frac{tp_s}{1-p_0-p_r} + i\left(\frac{p_0}{p_0+p_r} - \frac{p_s}{1 - p_0 -p_r}\right), i = 0, 1, ..., t.$$
$$\mathbb{E}(M_i(t)) = t(p_i + p_0), \quad  \,\mathbb{D}(M_i(t)) =t(p_i + p_0)(1 - p_i - p_0),\quad \frac{\mathbb{D}(M_i(t))}{\mathbb{E}(M_i(t))} = 1 - p_i - p_0, \quad   \,i = 1, 2,$$
$$\mathbb{E}(M_1(t)M_2(t)) = t(t-1)(p_1 + p_0)(p_2 + p_0) + tp_0,$$
 $$cov(M_1(t),M_2(t)) = t[p_0 - (p_1 + p_0)(p_2 + p_0)]$$

 Their correlation does not depend on $t$,
 $$cor(M_1(t),M_2(t)) = \frac{p_0 - (p_1 + p_0)(p_2 + p_0)}{\sqrt{[p_1 + p_0 - (p_1 + p_0)^2][p_2 + p_0 - (p_2 + p_0)^2]}}.$$

In general the processes $M_1$ and $M_2$ are not independent. However, for $p_0 = (p_1 + p_0)(p_2 + p_0)$, and $p_0 + p_1 \in [0, 1]$, which is possible to happen for example when $p_0 = p_1 = p_2 = \frac{1}{4}$,  the process $(M_1, M_2)$ reduces to a bivariate counting process which coordinates are independent binomial processes with parameters correspondingly $p_0 + p_1$ and $p_0 + p_2$.

Note also the formula for correlation  is not a particular case of formula (4.4) in Eagleson (1964) \cite{Eagleson1964},  as far as in our case the processes $B_1(t), B_2(t)$ and $B_0(t)$ are always dependent.

{\it Note  2.} The process $A = M_1 + M_2$ is a compound binomial process. For all $t = 0, 1, ...$, its time intersections satisfy the equalities
$$A(t) := M_1(t) + M_2(t) = \sum_{i = 0}^{B(t)} (I_{1i} + I_{2i} + 2I_{0i}) = \sum_{i = 0}^{t} (I_{A_i}(I_{1i} + I_{2i} + 2I_{0i})),$$
 where for any fixed $i = 1, 2, ...$, $\mathbb{P}(A_i) = p$, and $I_{A_i}$ are independent on $(I_{1i}, I_{2i}, I_{0i})$. It is easy to see that the process $A$ possesses the following properties. It is determined via the distribution of its time intersections which have the following p.g.fs.,
$$\mathbb{E}(z^{A(t)}) = \left[1 - p + (p_1 + p_2)z + p_0z^2\right]^t, \quad z \geq 0.$$
H.S.Steyn (1963) \cite{Steyn1963} calls these distributions "multivariate multinomial distributions". J.Panaretos and E. Xekalaki (1986) \cite{PanaretosXekalaki1986} use the name "cluster binomial distribution" and find its p.m.f. It is easy to see  that
$$\mathbb{E}(A(t)) = t(p_1 + p_2 + 2p_0), \quad \mathbb{D}(A(t)) = t[p_1 + p_2 + 4p_0 - t(p_1 + p_2 + 2p_0)^2].$$

In general $A$ is not a particular case of convolution of two binomial processes. However, for $p_0 = (p_1 + p_0)(p_2 + p_0)$, and $p_0 + p_1 \in (0, 1)$, which is possible to happen for example when $p_0 = p_1 = p_2 = \frac{1}{4}$ the process $A$ is sum of two independent binomial processes with parameters correspondingly $p_0 + p_1$ and $p_0 + p_2$.

\bigskip
\section{STOCHASTICALLY EQUIVALENT PRESENTATIONS OF THE TOTAL CLAIM AMOUNT PROCESS}

In this section we show that the process $(S_1, S_2)$ has random walk and compound binomial of type I, according to Sundt and Vernic (2009) \cite{SundtVernic}, presentations.

{\bf Theorem 2.} The bivariate total claim amount process $(S_1, S_2)$, described in (\ref{S1S2}) is a bivariate compound binomial process of type I. It is stochastically equivalent to the process $(S_3, S_4)$, where
\begin{eqnarray}
S_3(t) &=& \sum_{n = 0}^{B(t)}(I_{1n}Y_{1n} + I_{0n}Y_{3n}) = \sum_{n = 0}^{t}(I_{A_n}(I_{1n}Y_{1n} + I_{0n}Y_{3n})), \label{S3}\\
S_4(t) &=& \sum_{n = 0}^{B(t)}(I_{2n}Y_{2n} + I_{0n}Y_{4n}) = \sum_{n = 0}^{t}(I_{A_n}(I_{1n}Y_{2n} + I_{0n}Y_{4n})), \quad t = 0, 1, ....\label{S4}
\end{eqnarray}
Here $Y_{k0} = 0$, $k = 1, 2, 3, 4$, $(I_{11}, I_{21}, I_{01})$ is a vector with distribution (\ref{pmf1}), $B$ is a homogeneous binomial process with independent additive increments and parameter $p$. $I_{A_i}$, $i = 1, 2, ...$ are i.i.d. Bernoulli r.vs. with parameter $\mathbb{P}(A_i) = p$.  For any fixed $i = 1, 2, ...$, the random elements $I_{A_i}$, $B$, $Y_{1i}$, $(Y_{2i}, Y_{3i})$, $Y_{4i}$ and $(I_{1i}, I_{2i}, I_{0i})$ are mutually independent.

The equivalence is in the sense of their finite dimensional distributions.

{\bf Proof:} All processes have homogeneous and independent additive increments and they start from the coordinate beginning. Therefore, in order to prove their stochastic equivalence it is enough to prove equality in distribution of their univariate time intersections. Due to the uniqueness of the correspondence between the probability laws and their Laplace-Stieltjes transforms (LSTs), it is enough to derive equality between the LSTs.

Consider $t \in \mathbb{N}$. The definition of LSTs., (\ref{S1S2}),  the double expectation formula, applied with the distribution of $(B_1(t), B_2(t), B_0(t))$,  the multiplicative property of LSTs  and the definition for Multinomial distribution of $(B_1(t), B_2(t), B_0(t))$, entail
\begin{eqnarray*}
 && \mathbb{E}e^{-z_1S_1(t) - z_2S_2(t)} = \mathbb{E}e^{-z_1\sum_{i = 0}^{B_1(t)} Y_{1i} - z_1 \sum_{i = 0}^{B_0(t)} Y_{3i} - z_2 \sum_{i = 0}^{B_0(t)} Y_{4i}  - z_2 \sum_{i = 0}^{B_2(t)} Y_{2i}}\\
&=& \sum_{\substack{b_0, b_1, b_2 \in \{0, 1, ..., t\},\\ b_1 + b_2  + b_0 \leq t}} \mathbb{E}e^{-z_1 \sum_{i = 0}^{b_1} Y_{1i}}\mathbb{E}e^{-(z_1\sum_{i = 0}^{b_0} Y_{3i} + z_2\sum_{i = 0}^{b_0} Y_{4i})}\mathbb{E}e^{- z_2 \sum_{i = 0}^{b_2} Y_{2i}}\mathbb{P}(B_0(t)= b_0, B_1(t) = b_1, B_2(t) = b_2)\\
&=& \sum_{\substack{b_0, b_1, b_2 \in \{0, 1, ..., t\},\\ b_1 + b_2  + b_0 \leq t}} (\mathbb{E}e^{-z_1Y_{1i}})^{b_1}(\mathbb{E}e^{-(z_1Y_{3i} + z_2Y_{4i})})^{b_0}(\mathbb{E}e^{- z_2Y_{2i}})^{b_2}\mathbb{P}(B_0(t)= b_0, B_1(t) = b_1, B_2(t) = b_2)\\
&=& (1-p)^t\\
 &\times&\sum_{\substack{b_0, b_1, b_2 \in \{0, 1, ..., t\},\\ b_1 + b_2  + b_0 \leq t}}\left[\frac{p_1}{1-p}\mathbb{E}e^{-z_1Y_{1i}}\right]^{b_1}\left[\frac{p_0}{1-p}\mathbb{E}e^{-(z_1Y_{3i} + z_2Y_{4i})}\right]^{b_0}\left[\frac{p_2}{1-p}\mathbb{E}e^{- z_2Y_{2i}}\right]^{b_2} \frac{t!}{b_0!b_1!b_2!(t-b_0-b_1-b_2)!}  \\
&=& (1-p)^t \left[1 + \frac{p_1}{1-p}\mathbb{E}e^{-z_1Y_{11}}+  \frac{p_0}{1-p}\mathbb{E}e^{-(z_1Y_{31}+z_2Y_{41})}+\frac{p_2}{1-p}\mathbb{E}e^{-z_2Y_{21}}\right]^t\\
&=& \left[1 - p + p_1\mathbb{E}e^{-z_1Y_{11}} + p_0\mathbb{E}e^{-(z_1Y_{31}+z_2Y_{41})} + p_2\mathbb{E}e^{-z_2Y_{21}}\right]^t.
\end{eqnarray*}

From the other hand, the definitions (\ref{S3}), (\ref{S4}), the formula for double expectations, the independence between $B$ and the other components of the processes $S_3$ and $S_4$ entail
\begin{eqnarray*}
\mathbb{E}e^{-z_3S_3(t) - z_4S_4(t)} &=& \mathbb{E}e^{-z_3\sum_{i = 0}^{B(t)}(I_{1i}Y_{1i} + I_{0i}Y_{3i})-z_4\sum_{i = 0}^{B(t)}(I_{2i}Y_{2i} + I_{0i}Y_{4i})}\\
&=& \sum_{n = 0}^\infty \left[\mathbb{E}e^{-z_3I_{11}Y_{11} - (z_3Y_{31} + z_4Y_{41})I_{01}-z_4I_{21}Y_{21}}\right]^nP(B(t) = n).
\end{eqnarray*}

Now, (\ref{pmf1}), and the formula for p.g.f. of Binomial distribution entail
\begin{eqnarray*}
\mathbb{E}e^{-z_3S_3(t) - z_4S_4(t)} &=& \sum_{n = 0}^\infty \left[\frac{p_1}{p} \mathbb{E}e^{-z_3Y_{11}} + \frac{p_0}{p}\mathbb{E}e^{-(z_3Y_{31} + z_4Y_{41})} + \frac{p_2}{p} \mathbb{E}e^{-z_4Y_{21}}\right]^nP(B(t) = n)\\
&=& \left\{1 - p + p\left[\frac{p_1}{p} \mathbb{E}e^{-z_3Y_{11}} + \frac{p_0}{p}\mathbb{E}e^{-(z_3Y_{31} + z_4Y_{41})} + \frac{p_2}{p} \mathbb{E}e^{-z_4Y_{21}}\right]\right\}^t\\
&=& \left[1 - p + p_1\mathbb{E}e^{-z_3Y_{11}} + p_0\mathbb{E}e^{-(z_3Y_{31} + z_4Y_{41})} + p_2 \mathbb{E}e^{-z_4Y_{21}}\right]^{t}.
\end{eqnarray*}
which is exactly $\mathbb{E}e^{-z_1S_1(t) - z_2S_2(t)}$. In this way the proof is completed.  \hfill $\Box$

{\bf Corollary 3.} For all $t \geq 0$ the distribution of the time intersections $(S_1(t), S_2(t))$, $t = 1, 2, ...$ of the bivariate process $(S_1, S_2)$ defined in (\ref{S1S2}) are determined via the following LSTs,
$$\mathbb{E}e^{-z_1S_1(t) - z_2S_2(t)} = \left[1 - p + p_1\mathbb{E}e^{-z_1Y_{11}} + p_0\mathbb{E}e^{-(z_1Y_{31} + z_2Y_{41})} + p_2 \mathbb{E}e^{-z_2Y_{21}}\right]^{t}.$$
\begin{itemize}
\item[i)] If  $\mathbb{E}Y_{i1} < \infty$, $i = 1, 2, 3, 4$, then, $\mathbb{E}(S_1(t)) = t(p_1\mathbb{E}Y_{11} + p_0\mathbb{E}Y_{31})$, $\mathbb{E}(S_2(t)) = t(p_2\mathbb{E}Y_{21} + p_0\mathbb{E}Y_{41});$
\item[ii)] If $\mathbb{E}(Y_{i1}^2) < \infty$, $i = 1, 2, 3, 4$, then,
$$\mathbb{D}(S_1(t)) = t\left[p_1\mathbb{E}(Y_{11}^2) + p_0\mathbb{E}(Y_{31}^2) - (p_1\mathbb{E}Y_{11} + p_0\mathbb{E}Y_{31})^2\right],$$
$$\mathbb{D}(S_2(t)) =  t\left[p_2\mathbb{E}(Y_{21}^2) + p_0\mathbb{E}(Y_{41}^2) - (p_2\mathbb{E}Y_{21} + p_0\mathbb{E}Y_{41})^2\right];$$
\item[iii)]
$\mathbb{E}(S_1(t)S_2(t)) = tp_0\mathbb{E}(Y_{31}Y_{41})+ t(t-1)(p_1\mathbb{E}Y_{11} + p_0\mathbb{E}Y_{31})(p_2\mathbb{E}Y_{21} + p_0\mathbb{E}Y_{41})$;
\item[iv)]
$cov(S_1(t), S_2(t)) =  t\left[p_0\mathbb{E}(Y_{31}Y_{41})- (p_1\mathbb{E}Y_{11} + p_0\mathbb{E}Y_{31})(p_2\mathbb{E}Y_{21} + p_0\mathbb{E}Y_{41})\right];$
\item[v)] The correlation $cor(S_1(t),S_2(t))$  does not depend on $t$. More precisely
$$cor(S_1(t),S_2(t)) = \frac{ p_0\mathbb{E}(Y_{31}Y_{41})- (p_1\mathbb{E}Y_{11} + p_0\mathbb{E}Y_{31})(p_2\mathbb{E}Y_{21} + p_0\mathbb{E}Y_{41})}{\sqrt{\left[p_1\mathbb{E}(Y_{11}^2) + p_0\mathbb{E}(Y_{31}^2) - (p_1\mathbb{E}Y_{11} + p_0\mathbb{E}Y_{31})^2\right]\left[p_2\mathbb{E}(Y_{21}^2) + p_0\mathbb{E}(Y_{41}^2) - (p_2\mathbb{E}Y_{21} + p_0\mathbb{E}Y_{41})^2\right]}}.$$
\end{itemize}

{\bf Corollary 4.} The time intersections
$$S(t) := \sum_{i = 0}^{B(t)} \left[I_{1i}Y_{1i} + I_{0i}(Y_{3i}+Y_{4i}) + I_{2i}Y_{2i}\right], \quad t \geq 0,$$
of the compound binomial process $S = \{S(t), \,t\geq 0\}$, possess the following properties:
\begin{itemize}
\item[i)] their distributions are determined via the following probability generating functions,
$$\mathbb{E}(z^{S(t)}) = \left[1 - p + p_1 \mathbb{E} e^{-zY_{11}} + p_0\mathbb{E}e^{-z(Y_{31}+Y_{41})} + p_2 \mathbb{E}e^{-zY_{21}}\right]^{t}.$$
\item[ii)] $\mathbb{E}(S(t)) = t[p_1\mathbb{E}Y_{11} + p_0(\mathbb{E}Y_{31} + \mathbb{E}Y_{41})+ p_2\mathbb{E}Y_{21}]$;
\item[iii)] $\mathbb{D}(S(t)) = t\left\{p_1\mathbb{E}(Y_{11}^2)+ p_2\mathbb{E}(Y_{21}^2) + p_0\mathbb{E}[(Y_{31}+ Y_{41})^2] - [p_1\mathbb{E}Y_{11} + p_0(\mathbb{E}Y_{31} + \mathbb{E}Y_{41}) + p_2\mathbb{E}Y_{21}]^2\right\}.$
\end{itemize}

{\bf Theorem 3.} Consider the counting  process $(M_1, M_2)$ defined in (\ref{M}). If for all $i = 1, 2, ...$,  $Y_{1i} \bot Y_{2i} \bot Y_{3i} \bot Y_{4i}$, then the bivariate process $(S_1, S_2)$, described in (\ref{S1S2}) is stochastically equivalent to the process $(S_5, S_6)$, where
\begin{equation}\label{S5}
S_5(t) = \sum_{n = 0}^{M_1(t)}(Y_{1n} + Y_{3n}), \quad S_6(t) = \sum_{n = 0}^{M_2(t)}(Y_{2n} + Y_{4n}), \quad t = 0, 1, ....
\end{equation}
and the random process $(M_1, M_2)$, is independent on $Y_{1i}$, $Y_{2i}$, $Y_{3i}$ and $Y_{4i}$, $i = 1, 2, ...$.

The equivalence is in the sense of their finite dimensional distributions.

{\bf Proof:} Consider $t \in \mathbb{N}$. The definition of LSTs., (\ref{S1S2}),  the double expectation formula, applied with the distribution of $(M_1(t), M_2(t))$,  and the multiplicative property of LSTs, entail
\begin{eqnarray*}
 && \mathbb{E}e^{-z_1S_1(t) - z_2S_2(t)} = \mathbb{E}e^{-z_1\sum_{i = 0}^{M_1(t)}(Y_{1i} + Y_{3i}) - z_2 \sum_{i = 0}^{M_2(t)}(Y_{2i} + Y_{4i})}\\
 &=& \sum_{m_1 = 0}^t \sum_{m_2 = 0}^{t-m_1} (\mathbb{E}e^{-z_1Y_{1i}})^{m_1} (\mathbb{E}e^{-z_1Y_{3i}})^{m_1} (\mathbb{E}e^{-z_1Y_{4i}})^{m_2}(\mathbb{E}e^{-z_2Y_{4i}})^{m_2} \mathbb{P}(M_1(t) = m_1, M_2(t) = m_3)\\
 &=&  \left[1 - p + p_1\mathbb{E}e^{-z_1Y_{11}} + p_0\mathbb{E}e^{-z_1Y_{31}}\mathbb{E}e^{-z_2Y_{41}} + p_2 \mathbb{E}e^{-z_2Y_{21}}\right]^{t}
 \end{eqnarray*}
 which is exactly $\mathbb{E}e^{-z_1S_1(t) - z_2S_2(t)}$, where we have used Corollary 3 for  $Y_{3i} \bot Y_{4i}$.  \hfill $\Box$

\bigskip
\section{STOCHASTICALLY EQUIVALENT RISK MODELS}
Now, we are ready to show that the risk process $R$, defined in (\ref{R}) is a classical compound Binomial risk process. Therefore, in order to characterize the probabilities for ruin, the time of ruin, the surplus immediately before ruin, and the deficit at ruin given that ruin occurs, we can use the well-known results about the Binomial risk model or the discrete time risk reserve process. The last theory together with the theory of random walk could be seen in many textbooks in random processes, for example Rolski et al. 1998 \cite{Rolski1998} or Feller (1957, 2008) \cite{Feller1957,Feller2008} among others. By using these results we obtain Theorem 4, below.

The definition of the risk process $R_u$, the definitions of $S$, $S_1$ and $S_2$, (see \ref{S}, and \ref{S1S2}) and Theorem 2 entail that the risk process, described in (\ref{R}) is stochastically equivalent to the process
\begin{equation}\label{Rt}
R_u(t) = u + t - \sum_{i = 0}^{B(t)} Y_i = u - \sum_{i = 0}^{t} (I_{A_i}Y_i - 1), \quad t = 0, 1, ....
\end{equation}
where $Y_0 = 0$, and $Y_i: = I_{1i}Y_{1i} + I_{0i}(Y_{3i} + Y_{4i}) + I_{2i}Y_{2i}$. Therefore, if the r.vs $Y_{11}, Y_{2i}, Y_{3i}, Y_{4i}$ have integer-valued distribution, concentrated on some subset of $\mathbb{N}$, then
\begin{eqnarray}
\label{probf} \mathbb{P}(I_{A_i}Y_i - 1 = -1) &=& 1 - p \\
\nonumber \mathbb{P}(I_{A_i}Y_i - 1 = 0) &=& p\left(\frac{p_1}{p}\mathbb{P}(Y_{11} = 1) + \frac{p_2}{p}\mathbb{P}(Y_{12} = 1)\right) = p_1\mathbb{P}(Y_{11} = 1) + p_2\mathbb{P}(Y_{12} = 1)\\
\nonumber \mathbb{P}(I_{A_i}Y_i - 1 = k) &=& p_1\mathbb{P}(Y_{11} = k+1) + p_0\mathbb{P}(Y_{31} + Y_{41} = k+1) + p_2\mathbb{P}(Y_{21} = k+1), \,\, k \in \mathbb{N}.
\end{eqnarray}

Further on we assume that $\mu : = \mathbb{E} Y_1 = \frac{p_1}{p}\mathbb{E}Y_{11} + \frac{p_2}{p}\mathbb{E}Y_{21} + \frac{p_0}{p}(\mathbb{E}Y_{31} + \mathbb{E}Y_{41}) < \infty$. Then,
$$\mathbb{E}(I_{A_i}Y_i - 1) = p\mu - 1,$$
$$\mathbb{E}[R_u(t)] = u + t[1 - (p_1\mathbb{E}Y_{11} + p_2\mathbb{E}Y_{21} + p_0(\mathbb{E}Y_{31} + \mathbb{E}Y_{41}))] = u + t - p\mu t.$$

Now, we obtain the same safety loading like in  Th. Rolski et al. (1998) \cite{Rolski1998}, Shiu (1989) \cite{Shiu1989} or H. Gerber (1988) \cite{Gerber1988} for the corresponding compound Binomial model. More precisely,
\begin{equation}\label{safetyloading}
\rho:= \lim_{t \to \infty} \frac{\mathbb{E}[R_u(t)]}{\mathbb{E}[S(t)]} = \frac{1}{p\mu}-1,
\end{equation}
and therefore, the net profit condition  is
\begin{equation}\label{npc}
\rho > 0 \iff 1 > p\mu.
\end{equation}
It guaranties that in a long horizon the mean income of the insurer is bigger than the mean expenditures. Moreover, if the net profit condition (\ref{npc}) is satisfied, then the probability for ruin in "infinite horizon" and initial capital $u$, converges to $0$, when $u$ increases unboundedly and $\psi(u) < 1$, for all $u \in \mathbb{N}$. Otherwise, this probability is equal to $1$ for any initial capital $u \geq 0$. In the next theorem we apply the results for the discrete time risk model and obtain the corresponding results for our model.

We additionally denote the associated random walk and its global maxima by
\begin{equation}\label{RW}
\tilde{S}(t) := \sum_{i = 0}^{t} (I_{A_i}Y_i - 1), \quad t = 0, 1, ..., \quad Y_0 = 0,
\end{equation}
\begin{equation}\label{maxRW}
M:= max\{0,\tilde{S}(1), \tilde{S}(2), ...\}
\end{equation}
Then, $\delta(u) = \mathbb{P}(M \leq u)$.

The following lemma can be useful in order to prove that the distribution, described in Theorem 4. iii) coincides with the distribution, described in Theorem 4. iv). Its proof follows by the equality of the corresponding p.g.fs.

{\bf Lemma 1.} Let $X_1, X_2, ...$ be strictly positive i.i.d. r.vs. and $I_{A_1}, I_{A_2}, ...$ be i.i.d. Bernoulli r.vs. with $\mathbb{P}(A_1) = \frac{1-cp_G}{c(1-p_G)}$, where $p_G \in (0, 1)$, and $c$ is an arbitrary constant in the interval $(1, \frac{1}{p_G})$. Assume that
$$\mathbb{P}(\xi = k) = (1-cp_G)(cp_G)^k, \quad \mathbb{P}(\eta = k) = (1-p_G)p_G^k,\,\, k = 0, 1, ....$$
and all these r.vs. are independent. Then, 
$$\sum_{i=1}^{\xi} X_i \stackrel{d}{=} \sum_{i=1}^{\eta} (I_{A_i}X_i), \quad \sum_{i=1}^{0}:=0.$$

{\it Note:} This lemma shows that for any compound geometric distribution with non-negative summands there exists uncountably many  stochastically equivalent and compound geometric presentations.

The main theorem in this work is as follows.

{\bf Theorem 4.}
\begin{itemize}
  \item[i)]   If the net profit condition (\ref{npc}) is satisfied, $\mathbb{E}(z^{Y_{11}}) < \infty$, $\mathbb{E}(z^{Y_{31} + Y_{41}}) < \infty$, and $\mathbb{E}(z^{Y_{21}}) < \infty$, then
  $$\mathbb{E}(z^M) = \frac{(1 - p\mu)(1-z)}{1 - p + p_1 \mathbb{E}(z^{Y_{11}}) + p_0 \mathbb{E}(z^{Y_{31} + Y_{41}}) + p_2 \mathbb{E}(z^{Y_{21}}) - z}, \quad |z| < 1;$$
 \item[ii)]  If the net profit condition (\ref{npc}) is satisfied, $\delta(0) = \frac{1 - p\mu}{1 - p};$
  \item[iii)]  If the net profit condition (\ref{npc}) is satisfied, the distribution of $M$ is compound geometric. The number of summands $N$ has a p.m.f.  
  $$\mathbb{P}(N = k) = (p\mu)^k(1-p\mu), \,\,k = 0, 1, ...,$$
  and it is equal to $0$ otherwise. 
   
  The summands $U$ has p.m.f. 
  $$\mathbb{P}(U = k) = \frac{p_1\mathbb{P}(Y_{11} > k) + p_0\mathbb{P}(Y_{31} + Y_{41} > k) + p_2 \mathbb{P}(Y_{21} > k)}{p^2\mu},\,\, k = 0, 1, ...$$
     \item[iv)]  If the net profit condition (\ref{npc}) is satisfied, the distribution of $M$ is compound geometric.
    The number of summands $N_1$ has a p.m.f.
  $$\mathbb{P}(N_1 = k) = \left(\frac{1-p\mu}{1-p}\right)^k\left(\frac{p(\mu-1)}{1-p}\right), \,\,k = 0, 1, ...,$$
  and it is equal to $0$ otherwise. 
  
  The summands $U_1$ has p.m.f. 
  \begin{equation}\label{iv}
  \mathbb{P}(U_1 = k) = \frac{\mathbb{P}(Y_1 > k)}{\mu-1} = \frac{p_1\mathbb{P}(Y_{11} > k) + p_0\mathbb{P}(Y_{31} + Y_{41} > k) + p_2 \mathbb{P}(Y_{21} > k)}{p(\mu-1)}, \,\, k = 1, 2, ...
  \end{equation}
  \item[v)] If the net profit condition (\ref{npc}) is satisfied,  
  $$\delta(u) = \frac{1-p\mu}{1-p}\sum_{n=0}^\infty\left(\frac{p(\mu-1)}{1-p}\right)^nH^{*n}(u),$$
  $$\psi(u) = \frac{1-p\mu}{1-p}\sum_{n=1}^\infty\left(\frac{p(\mu-1)}{1-p}\right)^n[1-H^{*n}(u)],$$
  where $H(u)$ is the c.d.f. of the p.m.f. (\ref{iv}).
     \item[vi)]  If the net profit condition (\ref{npc}) is satisfied, and there exists $z > 1$ such that $1 - p + p_1 \mathbb{E}(z^{Y_{11}}) + p_0 \mathbb{E}(z^{Y_{31} + Y_{41}}) + p_2 \mathbb{E}(z^{Y_{21}}) = z$, then $\{z^{\sum_{i=1}^n I_{A_i}Y_i - n}, \, n = 0, 1, ...\}$ is a martingale with respect to the natural filtration, and
      $$\psi(u) = \frac{1}{\mathbb{E}[z^{\sum_{i=1}^{\tau(u)} I_{A_i}Y_i - \tau(u)}|\tau(u)<\infty]}  \leq \frac{1}{z^u}, \quad \forall u \geq 0.$$
    Here $\varepsilon = log(z)$ is the corresponding Lundberg exponent.
  \item[vii)] Consider $\lambda_r = \mathbb{P}(-R_0(\tau(0)) = r|\tau(0)< \infty)$. Then, for $r \in \mathbb{N}$,
  \end{itemize}
$$\lambda_{r+1} =  \frac{\lambda_r}{1 - p}[1 - p_1\mathbb{P}(Y_{11} = 1) - p_2\mathbb{P}(Y_{12} = 1) ] - \lambda_r\lambda_1- \frac{p_1\mathbb{P}(Y_{11} = r+1) + p_0\mathbb{P}(Y_{31} + Y_{41} = r+1) + p_2\mathbb{P}(Y_{21} = r+1)}{1 - p},$$
and these probabilities are equal to $0$, otherwise.
\begin{itemize}
  \item[viii)] If the net profit condition (\ref{npc}) is satisfied, then
\begin{equation}\label{expectation}
\mathbb{E}(-R_0(\tau(0))|\tau(0)< \infty) = \frac{p\mu -1}{p - (1-p)\lambda_1  - p_1\mathbb{P}(Y_{11} = 1) - p_2\mathbb{P}(Y_{12} = 1) }.
\end{equation}
\begin{equation*}
\mathbb{E}(\tau(0)|\tau(0)< \infty) = \frac{1}{p - (1-p)\lambda_1  - p_1\mathbb{P}(Y_{11} = 1) - p_2\mathbb{P}(Y_{12} = 1)}.
\end{equation*}
\item[iv)] If the net profit condition (\ref{npc}) is satisfied, then
$$\sum_{r=1}^\infty s^r\lambda_r = \frac{\sum_{r=1}^\infty s^r[p_1\mathbb{P}(Y_{11} = r+1) + p_0\mathbb{P}(Y_{31} + Y_{41} = r+1) + p_2\mathbb{P}(Y_{21} = r+1)] - \lambda_1q}{1 - p_1\mathbb{P}(Y_{11} = 1) - p_2\mathbb{P}(Y_{12} = 1) - \lambda_1q - \frac{q}{s}}.$$
  \item[x)] ({\bf{Distribution of the deficit at ruin, when $u = 0$}}) If the net profit condition (\ref{npc}) is satisfied, $p_1\mathbb{P}(Y_{11} = 1) + p_2\mathbb{P}(Y_{12} = 1) \geq 2p - 1$, and $\lambda_1 = \frac{p -  p_1\mathbb{P}(Y_{11} = 1) - p_2\mathbb{P}(Y_{12} = 1)}{1 - p}$, then
$$\lambda_r = \frac{\mathbb{P}(I_{A_1}Y_1 - 1 \geq r)}{1 - p} = \frac{p\mathbb{P}(Y_1 > r)}{1 - p}.$$
\end{itemize}

{\bf Proof:} {\bf i)} is an immediate corollary of Theorem 5.1.1. of Rolski et al. (1998) \cite{Rolski1998}, applied to the compound Binomial presentation of this model.

{\bf ii)} follows by (2.14) in Shiu (1989) \cite{Shiu1989} or by Frosting \cite{Frosting}, or by i) and the equality
 $$\delta(0) = \mathbb{P}(M \leq 0) = \mathbb{P}(M = 0) = \lim_{z \searrow 0} \mathbb{E}(z^M) =  \frac{1 - p\mu}{1 - p}.$$

{\bf iii)} follows by Theorem 5.1.1.(b) in Rolski et al. (1999) \cite{Rolski1998}.

{\bf iv)} is a corollary of the corresponding result in the compound Binomial risk model or we find the p.g.f. of this distribution, and we observe that it is exactly the one, described in i). The rest follows by the uniqueness of the correspondence between p.m.f. and p.g.f.

{\bf v)} follow by Shiu(1989) \cite{Shiu1989}.

{\bf vi)} Let  $\varepsilon = log(z)$. From the exponential change of measure (see e.g. Cor 3.5. p. 70 in Asmussen and Albrecher (2010) \cite{AsmussenAlbrecher2010}) if $z > 1$ is such that $\mathbb{E}z^{I_{A_1}Y_1 - 1} = 1$, then $\{z^{\sum_{i=1}^n I_{A_i}Y_i - n}, \, n = 0, 1, ...\}$ is a martingale with respect to the natural filtration. Now, Pr. II.3.1. and Cor. II.3.4 in Asmussen and Albrecher (2010) \cite{AsmussenAlbrecher2010}), and the equality  +
$$\mathbb{E}z^{I_{A_1}Y_1 - 1} = 1 - p + p_1 \mathbb{E}(z^{Y_{11}}) + p_0 \mathbb{E}(z^{Y_{31} + Y_{41}}) + p_2 \mathbb{E}(z^{Y_{21}}) - z$$
entail that
$$\psi(u) = \frac{z^{-u}}{\mathbb{E}[z^{\sum_{i=1}^{\tau(u)} I_{A_i}Y_i - \tau(u) - u}|\tau(u)<\infty]},$$
and $\psi(u) \leq e^{-\varepsilon u} = z^{-u}$ , for all $u \geq 0$.

{\bf vii)} Analogously to formula (3), in Frosting \cite{Frosting} we use the total probability formula with the hypothesis $"I_{A_1}Y_1 - 1 = -1"$, $"I_{A_1}Y_1 - 1 = 0"$, $"I_{A_1}Y_1 - 1 = k"$, $k \in \mathbb{N}$ and the Markov property of the random walks and obtain
\begin{equation}\label{lambda}
\lambda_r =  \mathbb{P}(I_{A_1}Y_1 - 1 = -1)(\lambda_{r+1} + \lambda_1\lambda_r) +  \mathbb{P}(I_{A_1}Y_1 - 1 = 0)\lambda_r  +\mathbb{P}(I_{A_1}Y_1 - 1 = r) \end{equation}
$$\lambda_r =  p_1\mathbb{P}(Y_{11} = r+1) + p_0\mathbb{P}(Y_{31} + Y_{41} = r+1) + p_2\mathbb{P}(Y_{21} = r+1) + [p_1\mathbb{P}(Y_{11} = 1) + p_2\mathbb{P}(Y_{12} = 1)]\lambda_r + (1 - p)(\lambda_{r+1} + \lambda_1\lambda_r).$$
Now we express $\lambda_{r+1}$ and complete the proof of v).

{\bf viii)} We multiply both sides of (\ref{lambda}) by $r+1$, then we sum up in both sides of this equality for $r \in \mathbb{N}$ and obtain (\ref{expectation}). The second expression follows by the fact that $-R_0(\tau(0)+)$ is a randomly stopper random walk; the facts that the summand have finite expectation;  $\tau(0)$ is a Markov time with respect to the natural filtration born by $\{\sum_{i=1}^n(I_{A_i}Y_i - 1), \,\, i = 0, 1, ...\}$;  given $\tau(0) < \infty$ this stopping time is almost sure finite, and the Wald equality. Finally, we use that  $\mathbb{E}(I_{A_i}Y_i - 1) = p\mu-1 $ and complete this part of the proof.

{\bf ix)} follows by v) when we multiply both sides of (\ref{lambda}) by $z^{r+1}$, and sum up them for $r \in \mathbb{N}$.

{\bf x)} is a corollary of v) and a particular case of Proposition 3, in Frosting \cite{Frosting}.

\hfill $\Box$

\bigskip
\section{CONCLUSIONS}

Having in mind Lemma 1 and the fact that the global maxima of a random walk is a compound geometrically distributed in this lemma and in Theorem 4 we show that:
\begin{description}
\item[$\diamond$] Although the p.g.f. uniquely determines the p.m.f. and vice versa, there exist uncountably many compound geometric and stochastically equivalent presentations of any compound geometric distribution with non-negative summands.
\item[$\diamond$] The global maxima of a random walk has uncountably many stochastically equivalent compound geometric presentations.
\item[$\diamond$] Probability to survive in much more general settings that those, discussed here, for example the Anderson risk model, has uncountably many Beekman's convolution series presentations.
\end{description}

Many analogous results for the random variables described in this risk process can be obtained as corollaries of the corresponding results in the classical compound Binomial risk model. One can see for example, Li (2005) \cite{Li2005}, for the distribution of the surplus before ruin and the one of the claim causing ruin, or Willmot (1993) \cite{Willmoth1993}, or Li and Sendova (2013) \cite{LiandSendova2013} for probabilities of ruin in finite time. The joint distribution and numerical characteristics of these random variables could be seen e.g. in Li et al. (2013) \cite{Lietal2013}.

\bigskip
\section{ACKNOWLEDGMENTS}
The first author is grateful to the Project RD-08-75/27.01.2021 from the Scientific Research Fund in Konstantin Preslavsky University of Shumen, Bulgaria.
The second author thanks to project No 2021 - FNSE – 05, financed by the Scientific Research Fund of Ruse University.
\nocite{*}
\bibliographystyle{aipnum-cp}%

\end{document}